\documentstyle[12pt]{article}

\def\be{\begin{equation}}
\def\ee{\end{equation}}
\def\IP{\hbox{\rm I\kern -1.6pt{\rm P}}}
\def\IC{{\hbox{\rm C\kern-.58em{\raise.53ex\hbox{$\scriptscriptstyle|$}} 
    \kern-.55em{\raise.53ex\hbox{$\scriptscriptstyle|$}} }}}
\def\IN{\hbox{I\kern-.2em\hbox{N}}}
\def\IR{\hbox{\rm I\kern-.2em\hbox{\rm R}}}
\def\ZZ{\hbox{{\rm Z}\kern-.33em{\rm Z}}}
\def\IT{\hbox{\rm T\kern-.38em{\raise.415ex\hbox{$\scriptstyle|$}} }}
\def\ang{\hbox{$<$\kern-.42em\hbox{\rm )} }}

\begin{document}

\title{Measures with infinite Lyapunov exponents\\ 
for the periodic Lorentz gas}
\author{N.I. Chernov
\\ Department of Mathematics\\
University of Alabama at Birmingham\\
Birmingham, AL 35294\\
chernov@vorteb.math.uab.edu
\and S. Troubetzkoy
\\ Department of Mathematics\\
State University of New York at Stony Brook\\
Stony Brook, NY 11794\\
serge@math.sunysb.edu\\
}
\date{\today}
\maketitle

\newtheorem{theorem}{Theorem}
\newtheorem{lemma}[theorem]{Lemma}
\newtheorem{corollary}[theorem]{Corollary}

\begin{abstract}  
We study invariant measures for the Lorentz gas which are supported on
the set of points with infinite Lyapunov exponents.  We construct examples
of such measures which are measures of maximal entropy and ones which are
not.

\medskip\noindent
Keywords:  billiard, Lorentz gas, Lyapunov exponent, measure of maximal entropy.
\end{abstract}

\section{Introduction}

For a given dynamical system a measure of maximal entropy captures
information about the ``most chaotic''
part of the dynamics.  More precisely a measure of maximal entropy
is a probability measure invariant under the dynamics
whose metric entropy is equal to
the topological entropy of the given system.  For systems with
finite topological entropy much is known about measures of maximal entropy.
Axiom A diffeomorphisms  have a unique measure of maximal entropy on each
topologically transitive component and this measure is Markovian
\cite{Ma}. Recently in \cite{KT} it was shown that $C^{1+\alpha}$ 
diffeomorphisms of compact Riemannian manifolds have at most countably many 
ergodic measures of maximal entropy and these measures are Markovian.  
For systems with infinite topological entropy nothing in this direction 
has been known up to now, although the main tool used to derive 
the above results, Markov partitions, has been developed for some time
\cite{BSC90}.

In \cite{Ch91a} it was shown that the billiard ball map for the periodic Lorentz
gas has infinite topological entropy.  In this article we study the set
of points with infinite Lyapunov exponents.  Using the cell structure developed
in \cite{BSC90,Ku} we construct an ergodic invariant probability measure
with infinite topological entropy
supported on this set. Since the topological entropy is infinite this  
is a measure of maximal entropy.  From
the construction it is clear that there many such measures can coexist
on a single component of topological transitivity.  
We also construct an ergodic invariant probability measure with finite entropy
which is supported on this set showing that infinite exponents do not
necessarily lead to infinite entropy.

\section{Cell structure and Lyapunov exponents}

We study a periodic Lorentz gas on a plane. For simplicity, we 
assume that every fundamental domain of this gas contains a single 
round scatterer. Let the fundamental domain be a unit square 
and the scatterer be the circle of radius $r>0$ centered at 
the origin. Thus we get a periodic array of circles of radius 
$r$ centered at sites of the integral lattice. A pointlike 
particle moves freely at unit speed between the circles and 
reflects elastically off them. The circles are immovable and rigid. 

We assume that $r<1/2$, so that the moving particle is not trapped 
between four neighboring circles. Then the particle can move 
freely without collisions indefinitely; such Lorentz gases 
are said to have {\it no horizon}. 

Since the gas is periodic we can consider a fundamental domain:
the configuration space $Q$ of this system is the unit torus 
$0\leq q_1,q_2 <1$ without the disc $q_1^2+q_2^2 \le r^2$ 
(mod 1). The phase space is $M=Q\times S^1$. This is a billiard 
system of Sinai type (with dispersing boundary). The billiard 
ball map $T$ is defined on the surface 
$$
    M_1 = \{ (q,v)\in M:\, q\in\partial Q\ {\rm and}\ 
        (v,n)\geq 0\}
$$
where $n$ is the inward unit normal vector to the boundary $\partial Q$ 
of the domain $Q$. The map $T$ is simply the first return map 
on the surface $M_1$, i.e. it sends the particle at a reflection 
point to its next reflection. We introduce the coordinates $(s,\varphi)$ 
on $M_1$, where $s$ is the arc length on the circle $\partial Q$ 
and $\varphi$ is the angle between the vector $v$ and the above 
normal vector $n$ to the circle, $0\leq s < 2\pi r$ and $-\pi/2 
\leq \varphi \leq \pi/2$. Since $s$ is a cyclic coordinate, 
$M_1$ is a cylinder. The map $T$ on $M_1$ preserves 
the measure $d\nu = c_{\nu}\cos\varphi ds\, d\varphi$, where 
$c_{\nu}=(2\pi r)^{-1}$ is the normalizing factor.

Sinai \cite{Si70} was first to study the properties of the map 
$T$ in detail. He proved the $T$ is hyperbolic, i.e. has nonzero 
Lyapunov exponents a.e., and constructed stable and unstable 
fibers at a.e. point $x\in M_1$. He also developed a proof of 
ergodicity of $T$, which was later improved in \cite{BS73}. 
In addition to ergodicity, the mixing and K-property of $T$ 
was established in \cite{Si70,BS73}, and its Bernoulli 
property was proved in \cite{GO}. Another proof of ergodicity, 
which worked for multidimensional Lorentz gas as well, was 
provided in \cite{SC87}. Sinai \cite{Si70} derived a formula 
for the Kolmogorov entropy of the map $T$, which was later 
reproved and studied in \cite{Ch91b}. In particular, it was 
shown in \cite{Ch91b} that the entropy $h(T)$ has the following 
asymptotics as $r\to 0$:
$$
    h(T)=2\ln(1/r) + O(1).
$$
Markov partitions for the map $T$ were constructed in \cite{BSC90}. 
Those provide a symbolic representation of $T$ by topological 
Markov chains with countable alphabet. Based on Markov 
partitions, it was later shown in \cite{Ch91a} that the topological 
entropy of $T$ is infinite. In particular, the natural ergodic 
measure $\nu$ on $M_1$ is not a measure of maximal entropy, 
since its entropy is finite. The statistical properties of the map $T$ 
were studied in \cite{BSC91}: a stretched exponential bound 
on the decay of correlations was established and the central 
limit theorem along with its weak invariance principle was proved. 
Note that the periodic Lorentz gas with no horizon, apparently, 
displays a `superdiffusive' behavior, as was conjectured and 
explained in \cite{Bl}. 

We will use the `cell structure' of the surface $M_1$ described 
in detail in \cite{BSC90,Ku}. The map $T$ has a countable 
number of domains of continuity, which accumulate at a finite 
number of singular points, at which the time of the first return 
is infinite. We will call such points {\it supersingular}. 
For example, four points on the 
circle $\partial Q$ with coordinates $(0,\pm r)$ and $(\pm r,0)$ 
and with $\varphi=\pm \pi/2$ are such supersingular points (they 
produce eight supersingular points in $M_1$). There might be more 
supersingular points for small radius $r$. The domains of continuity 
of the map $T$ (we call them {\it cells}) 
form a fairly standard structure in the 
neighborhood of every supersingular point, independent of $r$ 
(the structure is the same if the scatterers are not 
necessarily circles but smooth convex domains on the 
torus). The structure of cells is shown in Fig.~1. 
We denote cells $A_n$, $n\geq 2$, where $n$ means that 
the first return time is about $n$ on the cell $A_n$. 
Fig.~2 shows which points are included in the cell $A_n$. 
The sizes of the cells are shown in Fig.~1. Here $O(n^{-a})$ 
means a value between $c^{-1}n^{-a}$ and $cn^{-a}$ for 
some $c>1$. Since there are only a finite number of 
supersingular points, we can assume that the value of $c$ is 
the same for all of them. 

The inverse map $T^{-1}$ also has a countable number of 
domains of continuity, which accumulate at the 
supersingular points. They have a symmetric form shown by 
dashed lines in Fig.~1. We denote them by $A_n'$, 
$n\geq 2$, and call `inverse cells'. Clearly, any 
cell $A_n$ is mapped by $T$ onto an inverse cell 
$A_n'$ with the same value of $n$ but located in 
the neighborhood of another supersingular point. Fig.~3 
shows how $A_n$ is mapped onto $A_n'$ under $T$. 

It is shown in \cite{BSC90,BSC91} that unstable directions 
for the map $T$ are continuous at every supersingular 
point and the limit $(ds^u,d\varphi^u)$ of the unstable directions  
is positive and finite: $0<d\varphi^u/ds^u<\infty$. Likewise, 
the limit of the stable directions at every supersingular point is 
negative and finite: $0>d\varphi^s/ds^s>-\infty$. 
Therefore, we have transversality of stable directions
and increasing sides of the reverse cell $A_n'$ and
transversality of unstable directions in the neighborhood
of supersingular point and decreasing (long) sides of the cells
$A_n$. Thus it is 
clear from Fig.~3 that for all points $x\in A_n$ the 
one-step expansion in the unstable direction has a factor 
$O(n^{3/2})$, and one-step contraction in the stable direction 
has a factor $O(n^{3/2})$. 

We will study points $x\in M_1$ such that $T^ix$ belongs to 
some cell $A_{n_i}$ (near some supersingular point) with large 
$n_i=n_i(x)>0$ for every $i\in\ZZ$. In other words, we will study 
points whose trajectories stay very close to supersingular points 
all the time. It is clear that the positive Lyapunov exponent 
of any such point is 
\begin{eqnarray}
    \Lambda_+(x)
      &=&\lim_{I\to\infty} \frac{3}{2I}\sum_{i=1}^I\ln n_i(x)+O(1)\nonumber\\
      &=&\lim_{I\to\infty} \frac{3}{2I}\sum_{i=1}^I\ln n_{-i}(x)+O(1) 
        \label{L}
\end{eqnarray}
In particular, if both limits in (\ref{L}) are infinite, 
then $\Lambda_+(x)=\infty$. Similar formulas hold for 
the negative Lyapunov exponent, $\Lambda_-(x)$. 

In virtue of (\ref{L}), any point $x\in M_1$ such that $n_i(x)\to 
\infty$ as $|i|\to\infty$ has infinite Lyapunov exponents. Such 
points form a Cantor-like nonempty set concentrated in the 
vicinity of supersingular points. Obviously, this set does not 
support any finite invariant measure, because its every trajectory 
is attracted by supersingular points. 

However, there are points $x\in M_1$ for which the sequence 
$\{n_i(x)\}$ has infinity as a limit point (both as $i\to\infty$ 
and $i\to -\infty$) but that sequence is `recurrent', i.e. 
for every $n\geq 2$ there are asymptotic frequencies 
\be
  p_n^{\pm}(x)=\lim_{I\to\infty} \frac{\#\{i\in [1,I]:\, n_{\pm i}(x)=n\} }{I} 
    \label{pn}
\ee
and $\sum_n p_n^+(x)=\sum_n p_n^-(x) = 1$. If the sequences 
$\{p_n^+(x)\}$ and $\{p_n^-(x)\}$ decay slowly enough, then 
the corresponding point $x$ will have infinite Lyapunov exponent. 
It is clear that the following condition is sufficient for infinite 
Lyapunov exponents:
\be
     \sum_n p_n^{\pm}(x)\ln n = \infty
       \label{pnn}
\ee

\section{Measures with infinite Lyapunov exponents}

Here we construct ergodic measures for the map $T$ with finite and 
infinite Kolmogorov entropy, such that a.e. point $x\in M_1$ 
with respect to those measures has infinite Lyapunov exponents. 

The construction starts with the following observation. There 
are positive constants $c>1$ and $n_\ast\geq 1$ such that every 
cell $A_n$ with $n>n_\ast$ intersects all inverse cells $A_m'$ 
near the same supersingular point with $c\sqrt{n}\leq m\leq c^{-2}n^2$ 
so that both longer sides of $A_n$ cross both longer sides of 
$A_m'$ (as shown in Fig.~3). This observation is based on 
Fig.~1 and 3 and was made in \cite{BSC90}. 

Let $n_0,n_1,n_2,\ldots$ be a sequence of integers such that 
$n_i>n_\ast$ for all $i\geq 0$ and $c\sqrt{n_i}\leq n_{i+1}\leq 
c^{-2}n_i^2$ for all $i\geq 0$. Then there is a sequence of 
cells $A_{n_i}$ such that $A_{n_i}\cap T^{-1}A_{n_{i+1}}\neq 
\emptyset$ for all $i\geq 0$, and the intersection
$$
    \cap_{i=0}^{\infty} T^{-i}A_{n_i}
$$
is a monotone curve in the cell $A_{n_0}$ which stretches 
from its top (short) side to its bottom (short) side. 

The same is true for inverse cells: there is a sequence 
of inverse cells $A_{n_i}'$ such that the intersection 
$$
    \cap_{i=0}^{\infty} T^iA_{n_i}'
$$
is a monotone curve stretching from the top short side 
of $A_{n_0}'$ to its bottom short side. 

Consequently, for any double-infinite sequence of integers 
$\{n_i\}$, $-\infty <i<\infty$, such that $n_i>n_\ast$ and 
$c\sqrt{n_i}\leq n_{i+1}\leq c^{-2}n_i^2$ for all 
$i\in\ZZ$ (note that this condition is symmetric: 
$c\sqrt{n_i}\leq n_{i-1}\leq c^{-2}n_i^2$) 
there is a sequence of cells $\{A_{n_i}\}$ 
such that the intersection 
\be
    \cap_{i=-\infty}^{\infty} T^{-i}A_{n_i}
      \label{cap}
\ee
is a single point in the cell $A_{n_0}$. 

We now fix an increasing sequence of integers, $N_0,N_1,N_2, 
\ldots$, such that $N_{i+1}=[c^{-4}N_i^2]+1$ for all 
$i\geq 0$ and $N_0\gg n_\ast$. It has the following two 
properties:

(i) for any double-infinite sequence $\{n_i\}$, 
$-\infty <i<\infty$ such that $n_i=N_{s(i)}$ 
with some $s(i)\in\ZZ^+$ 
and $|s(i)-s(i+1)|\leq 1$ for all $i\in\ZZ$. Note 
the intersection (\ref{cap}) is a single point 
in the cell $A_{n_0}$; 

(ii) for all $i\geq 0$ 
$$
     N_i\geq c^{-4(2^i-1)}N_0^{2^i}
$$
so that if $N_0$ is large enough, the sequence $\{N_i\}$ 
grows at the following super-exponential rate: 
\be
    N_i\geq 2^{2^i}
      \label{N22}
\ee

The set of double-infinite sequences $\{n_i\}$ described by the 
condition (i) above is a topological Markov chain 
with a countable number of states, which can be identified 
with $N_1,N_2,\ldots$. The allowed transitions from 
every state $N_i$ are the ones to $N_i$ itself and 
to the two neighboring states, $N_{i-1}$ and $N_{i+1}$. 
The only exception is the first state, $N_1$, from 
which the transitions to itself and to $N_2$ are allowed. 

We denote the collection of the above double-infinite 
sequences $\{n_i\}$ by $\Omega_1$. In virtue of the 
property (i) every sequence $\omega =\{n_i\}\in\Omega_1$ 
corresponds to a point $x=x(\omega)\in M_1$ 
defined by the intersection (\ref{cap}). The 
set of points 
$$
   \Omega_{1,M}=\{x(\omega):\, \omega\in\Omega_1\} 
$$
is a closed Cantor-like subset of $M_1$ invariant under 
$T$, i.e. $T\Omega_{1,M}=T^{-1}\Omega_{1,M}=\Omega_{1,M}$. 

Let $\mu_1$ be a Markov measure on the symbolic space 
$\Omega_1$ defined by the following transition probabilities: 
$\pi_{i+1,i}=1/3$ for all $i\geq 1$, $\pi_{i-1,i} 
=2/3$ for all $i\geq 2$ and $\pi_{1,1}=2/3$. (Here 
$\pi_{i,j}$ stands for the probability of transition 
from $N_j$ to $N_i$). This Markov measure is ergodic 
and mixing, its stationary distribution is $p(N_i) = 
1/2^i$ for $i\geq 1$. 

The measure $\mu_1$ projected from $\Omega_1$ down to $M$ 
generates an ergodic measure $\nu_1$ for the map $T$, 
which is concentrated on $\Omega_{1,M}$. By the ergodic 
theorem, for $\nu_1$-almost every point $x\in \Omega_{1,M}$ 
the asymptotic frequencies $p_n^{\pm}(x)$ defined by 
(\ref{pn}) exist and are equal to $p_n^{\pm}(x)=1/2^i$ if 
$n=N_i$ for some $i\geq 1$ and zero otherwise. It is then a 
simple calculation based on (\ref{pnn}) and (\ref{N22}) 
that the Lyapunov exponents are infinite a.e. in $M$ 
with respect to the ergodic measure $\nu_1$: 
$$
   \sum_n p_n^{\pm}(x)\ln n = \sum_i p_{N_i}^{\pm}\ln N_i 
     \geq \sum_i 2^{-i}\ln 2^{2^i} = \infty
$$

The measure $\nu_1$ constructed above has a finite entropy. 
Indeed, it is a Markov measure, and so its entropy, see, 
e.g., \cite{Ma}, is given by 
\be
    h = - \sum_i p_i\sum_j \pi_{ij}\log\pi_{ij}
       \label{hpip}
\ee
where $\pi_{ij}$ are the transition probabilities and $p_i$ 
is the stationary distribution. In fact, since for every 
state $i$ only two transition probabilities $\pi_{ij}$ 
are positive, as defined above, we have $h\leq\log 2$ 
for the measure $\mu_1$. 

The above construction of the Markov measure $\mu_1$ can 
be modified so that its entropy will be infinite and 
Lyapunov exponents will be still infinite a.e. We outline 
the construction below. 

We now consider all the double-infinite sequences $\{n_i\}$, 
$-\infty <i<\infty$, satisfying $n_i>n_\ast$ and $c\sqrt{n_i}\leq 
n_{i+1} \leq c^{-2}n_i$ for all $i\in\ZZ$, 
as defined above. We denote the set of 
these sequences by $\Omega_2$. Obviously, $\Omega_2$ is a 
topological Markov chain with a countable number of states 
which can be identified with $n_\ast+1, n_\ast+2,\ldots$. 
We will number these states by $1,2,\ldots$ so that the $i$th 
state is identified with 
$n_\ast+i$. Every sequence $\omega =\{n_i\}\in\Omega_2$ 
corresponds to a point $x=x(\omega)\in M_1$ defined by the 
intersection (\ref{cap}). The set of points 
$$
   \Omega_{2,M}=\{x(\omega):\, \omega\in\Omega_2\} 
$$
is a closed Cantor-like subset of $M_1$ invariant under 
$T$, i.e. $T\Omega_{2,M}=T^{-1}\Omega_{2,M}=\Omega_{2,M}$. 

We are going to find an ergodic and mixing Markov measure $\mu_2$ 
on $\Omega_2$ with transition probabilities $\pi_{ij}$ and with 
stationary distribution $p_i$ satisfying two conditions: 

(iii) its entropy given by (\ref{hpip}) is infinite;

(iv) one has
$$
    \sum_i p_i \ln(n_\ast+i) = \infty
$$
so that, by the condition (\ref{pnn}), the projection of 
the measure $\mu_2$ on $M_1$ will have infinite Lyapunov 
exponents a.e. 

The existence of Markov measures satisfying (iii) and (iv) 
is not based on the dynamics of the Lorentz gas or billiards, 
and we only sketch a proof. 
It is clear that (iv) is always satisfied if the probabilities 
$p_i$ decay slowly enough, for example, if $p_i>{\rm const} 
\cdot (i\ln^2i)^{-1}$. Next, we will take care of the condition (iii). 

Let a probability distribution, $||p_i||$, be given. We will 
show how to find a transition matrix, $\Pi=||\pi_{ij}||$, 
preserving the distribution $||p_i||$, 
so that for every state $i\geq 1$ all the positive  
transition probabilities $\pi_{ij} > 0$ satisfy 
\be 
     c\sqrt{n_\ast+i}\leq n_\ast+j \leq c^{-2}(n_\ast+i)^2
      \label{ni}
\ee
in which case the Markov measure defined by $||p_i||$ 
and $||\pi_{ij}||$ will be concentrated on $\Omega_2$. 
For any $k\geq 1$ denote 
$$
    q_k = p_k-p_{k+1}+p_{k+2}-p_{k+3}+\cdots
$$
We assume that the probabilities $p_i$ decrease monotonically, 
$p_1>p_2>\cdots$, and that 
\be
     \frac 12 p_k \le q_k \le \frac 34 p_k
       \label{pqk}
\ee
for every $k\geq 1$. To assure this, it is enough to assume that 
\be
     \frac 13 \leq \frac{p_{k+2}-p_{k+1}}{p_{k+1}-p_k} \leq 1 
          \label{pppp} 
\ee
for all $k\geq 1$, which simply means that $\{p_k\}$ decays without 
abrupt drops. 

   Now, for any $k\geq 1$ define a matrix 
of transition probabilities, $\Pi (k)=||\pi_{ij}(k)||$, by 
$\pi_{ii}=1$ for all $1\leq i <k$, $\pi_{kk}=q_k/p_k$, 
$\pi_{i,i+1}=q_{i+1}/p_i$ for all $i\geq k$ and $\pi_{i,i-1} 
=q_i/p_i$ for all $i\geq k+1$. It is easy to check that 
all $\Pi(k)$ preserve the same distribution $||p_i||$. 
Now, let $1=k_1 < k_2 < k_3 < \cdots$ be a sequence of 
numbers such that $k_m = [\bar{c} m]$ for some sufficiently 
large $\bar{c}>2$. We then define the matrix of transition 
probabilities $\Pi=||\pi_{ij}||$ by 
$$
   \Pi = \Pi(k_1) \cdot \Pi(k_2) \cdot \Pi(k_3) \cdots
$$
(here every entry $\pi_{ij}$ requires only a finite number of 
multiplications). If $\bar{c}$ and $n_\ast$ are large enough, this 
matrix clearly satisfies (\ref{ni}). Now, due to (\ref{pqk}) 
all the positive entries of the matrices $\Pi(k)$, $k\geq 1$, 
are not smaller than $1/4$. In that case, the entropy 
of the conditional distribution, 
$$
    h(i) = -\sum_j \pi_{ij}\log\pi_{ij}
$$
increases to infinity as $i\to\infty$. Moreover, it is 
bounded below by some increasing sequence $h(i)\geq H_i$, 
$H_i\to\infty$ as $i\to\infty$, independent on the 
stationary distribution $||p_i||$. (The sequence $H_i$ 
depends only on the value of $\bar{c}$ which is only 
determined by $c$ and $n_\ast$ in (\ref{ni}). Finally, we 
can find a probability distribution $||p_i||$ satisfying 
(\ref{pppp}) and such that $\sum_i p_iH_i=\infty$. (Obviously, 
such a distribution always exists, whatever the increasing 
sequence $\{H_i\}$). Hence the condition (iii) is provided. 

The Markov measure $\mu_2$ satisfying the conditions 
(iii) and (iv) projected from $\Omega_2$ down to $M_1$ 
induces a measure, $\nu_2$, concentrated on $\Omega_{2,M}$. 
A.e. point with respect to $\nu_2$ has infinite Lyapunov exponents. 
This measure has infinite entropy, and so it is a measure of 
maximal entropy. Note that for every mixing subshift of finite 
type with finite entropy \cite {Pa} the measure of maximal entropy 
is a Markov one, and its transition probabilities are positive 
for all (topologically) allowed transitions between states. 
The measure of maximal entropy that we have constructed here 
is also a Markov one, but some of its transition probabilities 
are zero even for topologically allowed transitions. This 
means that we did not use all of the available topological 
richness of the chain $\Omega_2$ defined above. The transitions 
that we have used were enough to make the entropy of the 
Markov measure infinite.

\begin{center}
Figure captions.
\end{center}
\begin{enumerate}

\item
Fig. 1. A supersingular point $S$ and a cell $A_n$ near it.

\item
Fig. 2. The outgoing vectors form the cell $A_n$.

\item
Fig. 3. The image of $A_n$ under the map $T$ is a reverse cell, $A_n'.$

\end{enumerate}

\end{document}